October 10, 2007

# The Levy-Gromov Isoperimetric Inequality in Convex Manifolds with Boundary

Frank Morgan
Department of Mathematics and Statistics
Williams College
Williamstown, Massachusetts 01267
Frank.Morgan@williams.edu

**Abstract**

We observe after Bayle and Rosales that the Levy-Gromov isoperimetric inequality generalizes to convex manifolds with boundary.

**1. Introduction.** The Levy-Gromov isoperimetric inequality ([G, §2.2], [BZ, 34.3.2]) provides a sharp lower bound on the perimeter required to enclose given volume in a closed manifold M in terms of a positive lower bound on the Ricci curvature, by comparison with the sphere. Our Theorem 1 shows that the result and its proof generalize to compact, convex manifolds with boundary. The heart of the proof is the observation that the region of prescribed volume is covered by rays normal to the enclosing, area-minimizing surface, yielding an inequality for the volume in terms of the surface and its mean curvature after Heintze and Karcher [HK]. In a convex manifold with boundary, the surface should meet the boundary of M orthogonally, the normal rays should still cover the region, and the same result should hold.

The case of smooth boundary was proved by Bayle and Rosales [BR, Thm. 4.8], using the second variation formula.

So let M be a compact, convex n-dimensional Riemannian manifold with boundary B, smooth on the interior. The (partially free boundary) isoperimetric problem seeks a region of prescribed volume fraction $0 < V < 1$ of least perimeter, not counting perimeter inside B. We review some results from geometric measure theory (see [M1], [M3]). Existence follows from standard compactness arguments. Let S be the closure of the part of the boundary of the region in the interior of M. At an interior point p, S has a tangent cone, area-minimizing without volume constraint; if such a cone is a plane, then S is a smooth, constant-mean-curvature hypersurface at p. (Actually it is known that the interior singular set has Hausdorff dimension at most n–8.) At boundary points of S, tangent "cones" still exist, and one expects S to be normal to B at most regular points of B, although we do not need such delicate regularity.

For our convexity hypothesis, all we need to know is that M has a convex tangent cone at all boundary points, that some shortest path between interior points lies on the interior, and that some shortest path between an interior point and a boundary point q is not tangent to the boundary at q.

**1. Theorem** (Levy-Gromov for convex manifolds with boundary). *Let M be an n-dimensional, compact, convex, connected Riemannian manifold with boundary, smooth on the interior, with Ricci curvature bounded below by n−1 and volume λ times the volume of the unit sphere. Then the isoperimetric profile P(V), the least perimeter to enclose given volume fraction 0<V<1, satisfies*

(9) $$P \geq \lambda P_0,$$

*where $P_0$ is the isoperimetric profile of the unit sphere. Equality holds only if M is the round suspension of a space form with boundary.*

We provide a sketch of the whole proof. The main new ingredient for the generalization appears in paragraph three, the observation that shortest paths from an interior point to an isoperimetric surface meet the surface on the interior, so that normal rays from the interior of the surface cover the interior of the region.

*Proof sketch.* For given volume fraction $0 < V < 1$, let P be the perimeter of a minimizing hypersurface S in M and let $P_0$ be the perimeter of the hypersphere $S_0$ in the unit sphere. By replacing V by 1−V (which changes the sign of the mean curvatures) if necessary, we may assume that the mean curvature of S is greater than or equal to that of $S_0$.

The idea is to estimate V by the volume of the union of rays normal to S at regular points of S. Consider a shortest path γ from a point p to the surface S, meeting S at a point q. We claim that if q is an interior point, then q must be a regular point of S. Since γ is a shortest path, an (area-minimizing) tangent cone C to S at q lies in the far halfspace bounded by the hyperplane normal to γ. It follows that C must be a hyperplane (or moving the vertex in the direction of the continuation of γ would decrease area to first order). Moreover, as a shortest path, γ must meet S normally.

We claim that q cannot be a boundary point; suppose it were. Since γ is a shortest path, C must lie in the far half-space intersect the convex tangent cone C' to the boundary of M at q. It follows that C must be the portion of a hyperplane inside C' and that the line tangent to γ at q must lie in C' (or moving the vertex in the direction of the continuation of γ would decrease area to first order; if the line leaves C' keep only the part of C inside C' and do better). Thus γ is tangent to the boundary of M at q, a contradiction of convexity.

Therefore the union of rays normal to S at regular points of S cover the region of volume fraction V. By the calculus estimate of Heintze-Karcher [HK], the volume λV enclosed by S satisfies

$$\frac{\lambda V}{P} \leq \frac{V}{P_0},$$

as desired. If equality holds for some V, then a pencil of geodesics normals to S must be isometric to a pencil of normals to the equator of a round sphere, M is the round suspension of a convex space form with boundary, and equality holds for all V.

*Technical note.* While the existence of a tangent cone depends on so-called monotonicity, we just need a non-zero tangent object (weak limit of homothetic expansions), which follows from lower density bounds (which can be proved for example as in [M2, Lemma 3.6]) and trivial upper density bounds. This not an issue when M is smooth up to the boundary.

**2. Corollary.** *Let R be a convex subregion of the nD unit hemisphere with volume λ times the volume of the hemisphere. Then the isoperimetric profile P(V) (least perimeter to enclose given volume fraction 0<V<1) satisfies*

(9)                              $P \geq \lambda P_1$ ,

*where $P_1$ is the isoperimetric profile of the unit hemisphere. If equality holds for some V, then M is the suspension of a convex subset of an equatorial hypersphere (and hence equality holds for all V).*

*Proof.* Corollary 2 follows immediate from Theorem 1 because $P_1 = .5 P_0$ and the λ of the corollary is twice the λ of the theorem.

As a further corollary we have an isoperimetric result stated in Morgan [M2, Rmk. 3.11], a paper which dealt primarily with the *surface* of polytopes. Corollary 3 was proved earlier without uniqueness by Lions and Pacella [LP, Thm. 1.1]. See also [M4, Rmk. after Thm. 10.6] and for the smooth case [RR,Thm. 4.11].

**3. Corollary**. *Let $P^n$ be a solid (convex) polytope in $\mathbf{R}^n$. For small prescribed volume, isoperimetric regions are balls about a vertex.*

*Proof.* The proof for solid polytopes follows the proof for surfaces of polytopes [M2]. The proof applies an isoperimetric comparison theorem [Ros, Thm. 3.7] for products and cones to the cone C over a vertex of P. Roughly, since balls about the origin are isoperimetric in a halfplane, they are isoperimetric in C. There is a hypothesis on the link, which requires Corollary 2 above, without the need for smoothing or the consequent loss of uniqueness which occur in the original proof for surfaces of polytopes.

# References

[BR]      Vincent Bayle and César Rosales, Some isoperimetric comparison theorems for convex bodies in Riemannian manifolds, Indiana Univ. Math. J. 54 (2005), 1371-1394.

[BZ]      Yu. D. Burago and V. A. Zalgaller, Geometric Inequalities, Springer-Verlag, 1988.

[G]     M. Gromov, Isoperimetric inequalities in Riemannian manifolds, Appendix I to Vitali D. Milman and Gideon Schechtman, Asymptotic Theory of Finite Dimensional Normed Spaces, Lec. Notes Math. 1200, Springer-Verlag, 1986.

[HK]     Ernst Heintze and Hermann Karcher, A general comparison theorem with applications to volume estimates for submanifolds, Ann. Scient. Éc. Norm. Sup. 11 (1978), 451-470.

[M1]     Frank Morgan, Geometric Measure Theory: a Beginner's Guide. Academic Press, third edition, 2000, fourth edition, to appear.

[M2]     Frank Morgan, In polytopes, small balls about some vertex minimize perimeter, J. Geom. Anal. 17 (2007), 97-106.

[M3]     Frank Morgan, Regularity of isoperimetric hypersurfaces in Riemannian manifolds, Trans. Amer. Math. Soc. 355 (2003), 5041-5052.

[M4]     Frank Morgan, Riemannian Geometry: a Beginner's Guide, A. K. Peters, Natick, Massachusetts, 1998.

[RR]     Manuel Ritoré and César Rosales, Existence and characterization of regions minimizing perimeter under a volume constraint inside Euclidean cones, Trans. Amer. Math. Soc. 356 (2004), 4601-4622.

[Ros]     Antonio Ros, The isoperimetric problem, David Hoffman, ed., Global Theory of Minimal Surfaces (Proc. Clay Math. Inst. 2001 Summer School, MSRI), Amer. Math. Soc., 2005, 175-209. (Preprint at http://www.ugr.es/~aros/isoper.pdf.)